# Theory on the Structure and Coloring of Maximal Planar Graphs

## (1)Recursion Formulae of Chromatic Polynomial and Four-Color Conjecture

XU Jin

(*School of Electronic Engineering and Computer Science*, *Peking University*, *Beijing* 100871, *China*)

(*Key Laboratory of High Confidence Software Technologies*, *Peking University*, *Beijing* 100871, *China*)

**Abstract**: In this paper, two recursion formulae of chromatic polynomial of a maximal planar graph $G$ are obtained: when $\delta(G)=4$, let $W_4^v$ be a 4-wheel of $G$ with wheel-center $v$ and wheel-cycle $v_1v_2v_3v_4v_1$, then $f(G,4)=f((G-v)\circ\{v_1,v_3\},4)+f((G-v)\circ\{v_2,v_4\},4)$; when $\delta(G)=5$, let $W_5^v$ be a 5-wheel of $G$ with wheel-center $v$ and wheel-cycle $v_1v_2v_3v_4v_5v_1$, then $f(G,4)=[f(G_1,4)-f(G_1\cup\{v_1v_4,v_1v_3\},4)]+[f(G_2,4)-f(G_2\cup\{v_3v_1,v_3v_5\},4)]+[f(G_3,4)-f(G_3\cup\{v_1v_4\},4)]$, $G_1=(G-v)\circ\{v_2,v_5\}$, $G_2=(G-v)\circ\{v_2,v_4\}$, $G_3=(G-v)\circ\{v_3,v_5\}$, where"$\circ$"denotes the operation of vertex contraction. Moreover, the application of the above formulae to the proof of Four-Color Conjecture is investigated. By using these formulae, the proof of Four-Color Conjecture boils down to the study on a special class of graphs, viz., 4-chromatic-funnel pseudo uniquely-4-colorable maximal planar graphs.

## 1 Introduction

All graphs considered in this paper are finite, simple and undirected. For a given graph $G$, we use $V(G)$, $E(G)$, $d_G(v)$ and $N_G(v)$ to denote the *vertex set*, the *edge set*, the *degree* of $v$ and the *neighborhood* of $v$ in $G$ (the set of neighbors of $v$), respectively, which can be written as $V$, $E$, $d(v)$ and $N(v)$ for short. The *order* of $G$ is the number of its vertices. A graph $H=(V',E')$ is a *subgraph* of $G$ if $V'\subseteq V$ and $E'\subseteq E$. For a subgraph $H$ of $G$, whenever $u,v\in V'$ are adjacent in $G$, they are also adjacent in $H$, then $H$ is an *induced subgraph* of $G$ or a subgraph of $G$ *induced* by $V'$, denoted by $G[V']$. Two graphs $G$ and $H$ are *disjoint* if they have no vertex in common. By starting with a disjoint union of $G$ and $H$, and adding edges joining every vertex of $G$ to every vertex of $H$, one obtains the *join* of $G$ and $H$, denoted $G\vee H$. We write $K_n$ and $C_n$ for the *complete graph* and *cycle* of order $n$, respectively. The join $C_n\vee K_1$ of a cycle and a single vertex is referred to as a *wheel*, denoted by $W_n$ (the examples $W_3,W_4,W_5$ are shown in Figure 1, where $C_n$ is called the *cycle* of $W_n$ and the vertex of $K_1$ is called the *center* of $W_n$. If $V(K_1)=\{x\}$, we also denote by $C^x$ the cycle of $W_n$. A graph is $k$*-regular* if all of its vertices have the same degree $k$. A 3-*regular* graph is usually called a *cubic graph*.

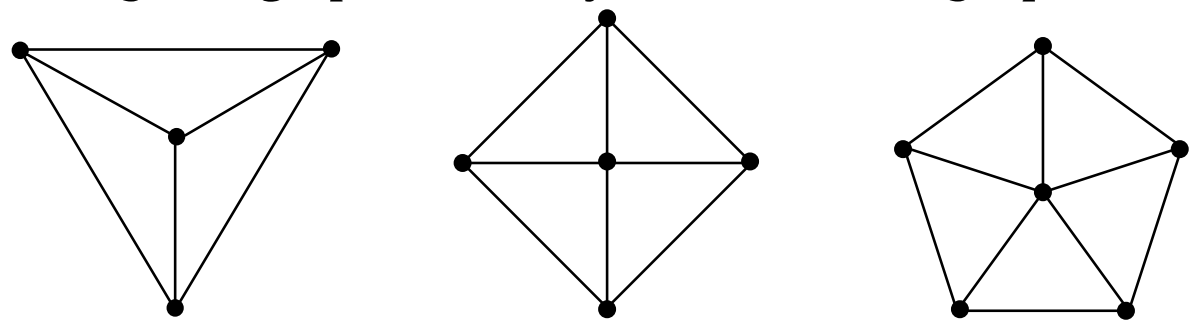

(a) $W_3$ (b) $W_4$ (c) $W_5$

Figure 1. Three wheels $W_3,W_4,W_5$.

A $k$*-coloring* of a graph $G$ is a mapping $f$ from the vertex set $V$ to the color set $C(k)=\{1,2,\cdots,k\}$ such that $f(x)\neq f(y)$ for any $xy\in E(G)$. A graph $G$ is $k$*-colorable* if it admits a $k$-coloring. The minimum $k$ for which a graph $G$ is $k$-colorable is called its *chromatic number*, denoted by $\chi(G)$. If $\chi(G)=k$, then $G$ is called a $k$-chromatic graph. Alternatively, each $k$-coloring $f$ of $G$ can be viewed as a partition $\{V_1,V_2,\cdots,V_k\}$ of $V$, where $V_i$ denotes the set of vertices assigned color $i$, called a *color class* of $f$. So it can be written as $f=(V_1,V_2,\cdots,V_k)$. In other words,

$$V(G)=\bigcup_{i=1}^{k}V_i,\ V_i\neq\varnothing,\ V_i\cap V_j=\varnothing,\ i\neq j,\ i,j=1,2,\cdots,k$$

where $V_i$ is an independent set of $G$, $i=1,2,\cdots,k$.

---

*Corresponding author: jxu@pku.edu.cn

Foundation Items: The National 973 Program of China (2013CB329600), The National Natural Science Foundation of China (61472012, 6152046, 6152012, 61572492, 61372191, 61472012)

The set of all $k$-colorings of a graph $G$ can be denoted by $C_k\ G$ . For a $k$-chromatic graph $G$, the notation $C_k^0\ G$ denotes the set of all the partitions of $k$-coloring class of $G$, simplified by the *partition set of $k$-color class* of $G$.

A graph is said to be *planar* if it can be drawn in the plane so that its edges intersect only at their ends. Such a drawing is called a *planar embedding* of the graph. Any planar graph considered in the paper is under its planar embedding. A *maximal planar graph* is a planar graph to which no new edges can be added without violating planarity. A *triangulation* is a planar graph in which every face is bounded by three edges (including its infinite face). It can be easily proved that maximal planar graphs are triangulations, and vice versa.

The planar graph is a very important class of graphs no matter which aspect, theoretical or practical, is concerned. In theory, there are many famous conjectures that have very significant effect on graph theory, even mathematics, such as the Four-Color Conjecture, the Uniquely Four-Colorable Planar Graphs Conjecture, the Nine-Color Conjecture and Three-Color Problem etc[1]. In application, planar graphs can directly be applied to the study of layout problems[2] and information science[3] etc.

Because the studying object of the well-known Four-Color Conjecture can be confined to maximal planar graphs, many scholars have been strongly attracted to the study of this typical topic. They did research on maximal planar graphs from a number of different standpoints, such as degree sequence, construction, coloring, traversability and generating operations, etc[4]. Moreover, many new conjectures on maximal planar graphs have been proposed, for instance, Uniquely Four-Colorable Planar Graphs Conjecture and Nine-Color Conjecture. These conjectures have gradually become the essential topics on maximal planar graphs.

In the process of studying Four-Color Conjecture, one important method, finding an unavoidable set of reducible configurations, was proposed. This method has been used in Kempe's "proof"[5], Heawood's counterexample[6] and the computer-assisted proof due to Appel and Haken[7-9]. Using this method, Appel and Haken found an unavoidable set containing 1936 reducible configurations and proved Four-Color Conjecture. In 1997, Robertson, Sanders, Seymour, Thomas, et al. [10,11] gave a simplified proof. They found an unavoidable set containing only 633 reducible configurations.

The research on unavoidable sets originated from Wernicke's work[12] in 1904. The concept of reducibility was introduced by Birkhoff[13] in 1913. On the research for finding an unavoidable set of reducible configurations, the great contribution was made by German mathematician Heesch[14]. He introduced a method "discharging" to find an unavoidable set of a maximal planar graph, which lied the foundation for solving Four-Color Conjecture by electronic computer in 1976[7-9]. Moreover, many researchers studied Four-Color Problem by this method, such as Franklin[15,16], Reynolds[17], Winn[18], Ore and Stemple[19] and Mayer[20].

However, these proofs were all computer-assisted and hard to be checked one by one by hand. Therefore, finding a mathematical method to concisely solve the Four-Color Problem is still an open hard problem.

Another incorrect proof of Four-Color problem[21] was given by Tait in 1880. His proof was based on an assumption: each 3-connected cubic plane graph was Hamiltonian. Because this assumption is incorrect, Tait's proof is incorrect. Although the error in his proof was found by Petersen[22] in 1898, the counterexample was not given until 1946[23]. Then, in 1968, Grinberg[24] obtained a necessary condition, thus producing many non-Hamiltonian cubic planar graphs of 3-connected. Although the proof of Tait was incorrect, his work had a strong influence on the research on Graph Theory, especially edge-coloring theory.

Let $f\ G,t$ be the number of colorings for the vertices of a labeled graph $G$ with $t$ colors. Obviously, if $t<\chi\ G$ , $G$ can not be properly colored, so $f\ G,t\ =0$. But if $\chi\ G\ \le t$, then $G$ admists this coloring must exist, and that is $f\ G,t$ $>0$. For every planar graph $G$, if $f\ G,4\ >0$ can be proved, it is equivalent to the proof of the Four-Color Problem! This is the method that Birkhoff[25,26] had proposed for attacking Four-Color Problem in 1912. Later on, it was found that $f(G,t)$ is a polynomial in the number $t$, called the *chromatic polynomial* of graphs, which has become a fascinating branch in the field of graph theory at present[27]. But it was a pity that Birkhoff's aim had not been reached. Further research on chromatic polynomials can be found in references [25–31]. The best result, due to Tutte[28], was that if $t=\tau(\sqrt{5})=3.618\cdots$ (where $\tau=(\sqrt{5}+1)/2$ ), then $f(G,\tau\sqrt{5})>0$. The result seemed to be a pity that it brushed past the Four-Color Problem, because the Four-Color Conjecture holds if $f(G,4)>0$.

In order to calculate the chromatic polynomial of a given graph, the basic tool is the *Deletion-Contract Edge Formula*.

For an edge $e$ of a graph $G$, we denote by $G-e$ and $G\circ e$ the graphs obtained from $G$ by deleting and contracting the edge $e$, respectively. Throughout the contraction operation, all graphs have no loops and parallel edges, except $W_2$.

**The Deletion-Contract Edge Formula.** For a given graph $G$ and an edge $e\in E(G)$, we have

$$f(G,t)=f(G-e,t)-f(G\circ e,t).$$

Moreover, the author[32,33] obtained a recursion formula of chromatic polynomial by vertex deletion and a chromatic polynomial between a graph and its complement.

Perhaps for the perfect degree of Tutte's work and his highly status in academia, once upon a time, it was thought that to attack the Four-Color Problem by chromatic polynomial is impossible. Nevertheless, our works below give a new hope to solve the Four-Color Problem by chromatic polynomial.

## 2 Recursion formulae of chromatic polynomial by contracting wheels

We first give two useful lemmas as follows.

**Lemma 1.**[26] For any planar graph $G$, it is 4-colorable if and only if

$$f(G,4)>0 \qquad (1).$$

**Lemma 2.**[25,27] Let $G$ be the union of two subgraphs $G_1$ and $G_2$, whose intersection is a complete graph of order $k$. Then

$$f(G,t)=\frac{f(G_1,t)\times f(G_2,t)}{t(t-1)\cdots(t-k+1)} \qquad (2)$$

**Theorem 1**. Let $G$ be a maximal planar graph, $v$ be a 4-degree vertex of $G$, and $N(v)=\{v_1,v_2,v_3,v_4\}$ (see Figure 2). Then

$$f(G,4)=f(G_1,4)+f(G_2,4) \qquad (3)$$

where $G_1=(G-v)\circ\{v_1,v_3\}$, $G_2=(G-v)\circ\{v_2,v_4\}$.

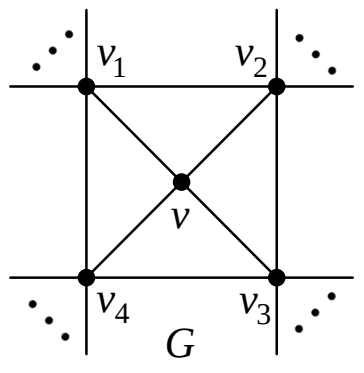

Figure 2. A maximal planar graph with a 4-degree vertex

**Proof.** In the following derivation, we represent $G$ by $G[\overline{N(v)}]$. Now we first compute the chromatic polynomial of the graph $G$ by the Deletion-Contract Edge Formula. For the sake of understanding clearly, a method introduced by Zykov[34] is used here, where the chromatic polynomials are represented by the corresponding graphical graphs without $t$. Notice that if there are at least two edges adjacent to two vertices, then only one remains and others are deleted excluding $W_2$.

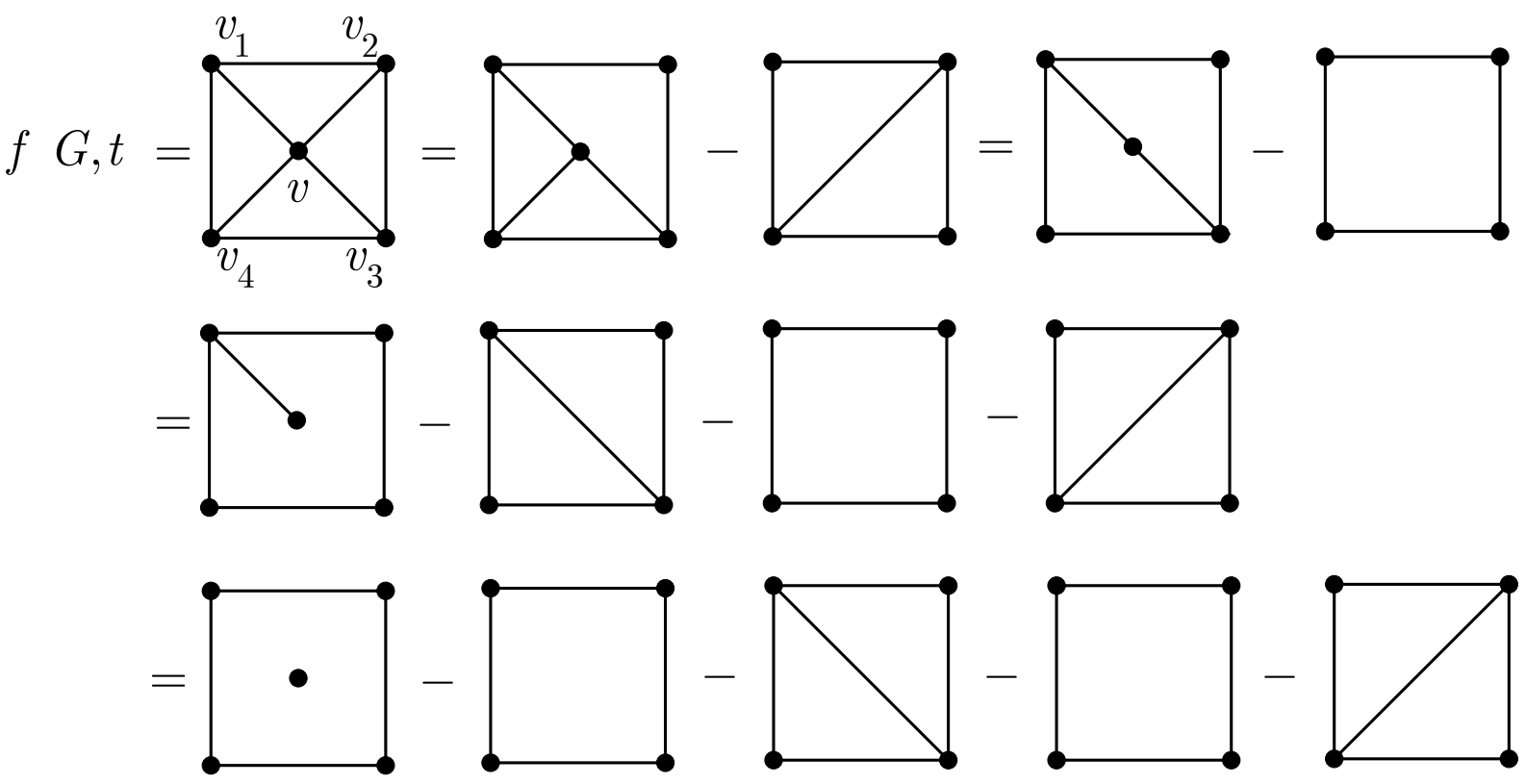

(4)

By Lemma 2, the chromatic polynomial of the first subgraph in Formula (4) is $t \cdot f(G-v,t)$. Therefore,

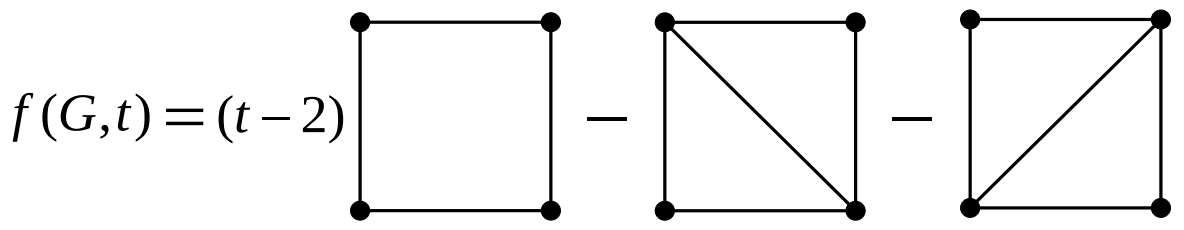

(5)

When $t=4$, we can obtain that

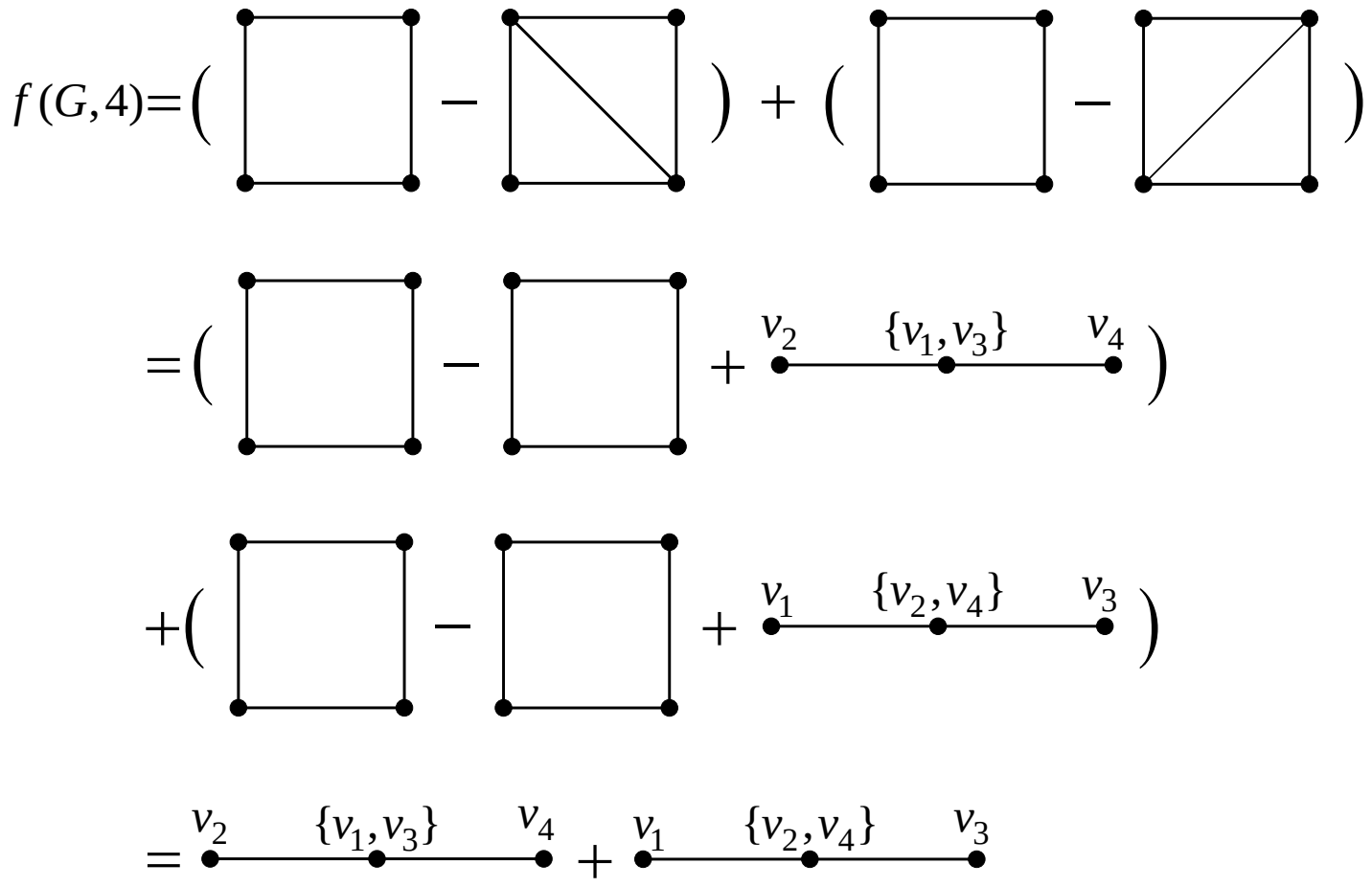

(6)

Notice that the two graphs in Formula (6) denote $(G-v)\circ\{v_1,v_3\}$ and $(G-v)\circ\{v_2,v_4\}$, respectively. It is easily proved that they are both maximal planar graphs of order $n-2$. Thus, we obtain that

$$f(G,4)=f((G-v)\circ\{v_1,v_3\},4)+f((G-v)\circ\{v_2,v_4\},4)=f(G_1,4)+f(G_2,4) \tag{7}$$

namely,

$$f(G,4)=f(G_1,4)+f(G_2,4) \tag{8}$$

**Theorem 2**. Let $G$ be a maximal planar graph, $v$ be a 5-degree vertex of $G$, and $N(v)=\{v_1,v_2,v_3,v_4,v_5\}$ (see Figure 3). Then

$$\begin{aligned} f(G,4) = & [f(G_1,4)-f(G_1\cup\{v_1v_4,v_1v_3\},4)] \\ & +[f(G_2,4)-f(G_2\cup\{v_3v_1,v_3v_5\},4)] \\ & +[f(G_3,4)-f(G_3\cup\{v_1v_4\},4)] \end{aligned} \tag{9}$$

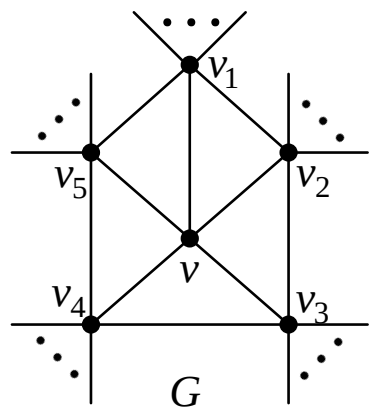

Figure 3. A maximal planar graph with a 5-degree vertex

where $G_1=(G-v)\circ\{v_2,v_5\}$, $G_2=(G-v)\circ\{v_2,v_4\}$, $G_3=(G-v)\circ\{v_3,v_5\}$.

**Proof.** The graph $G$ is represented by $G[\overline{N(v)}]$ in the following proof. The chromatic polynomial of graph $G$ can be calculated by applying the Deletion-Contract Edge Formula repeatedly. If parallel edges appear in the process, reserve only one edge excluding $W_2$. We use $W_5$ to represent the chromatic polynomial of $G$. In this way, we can obtain that

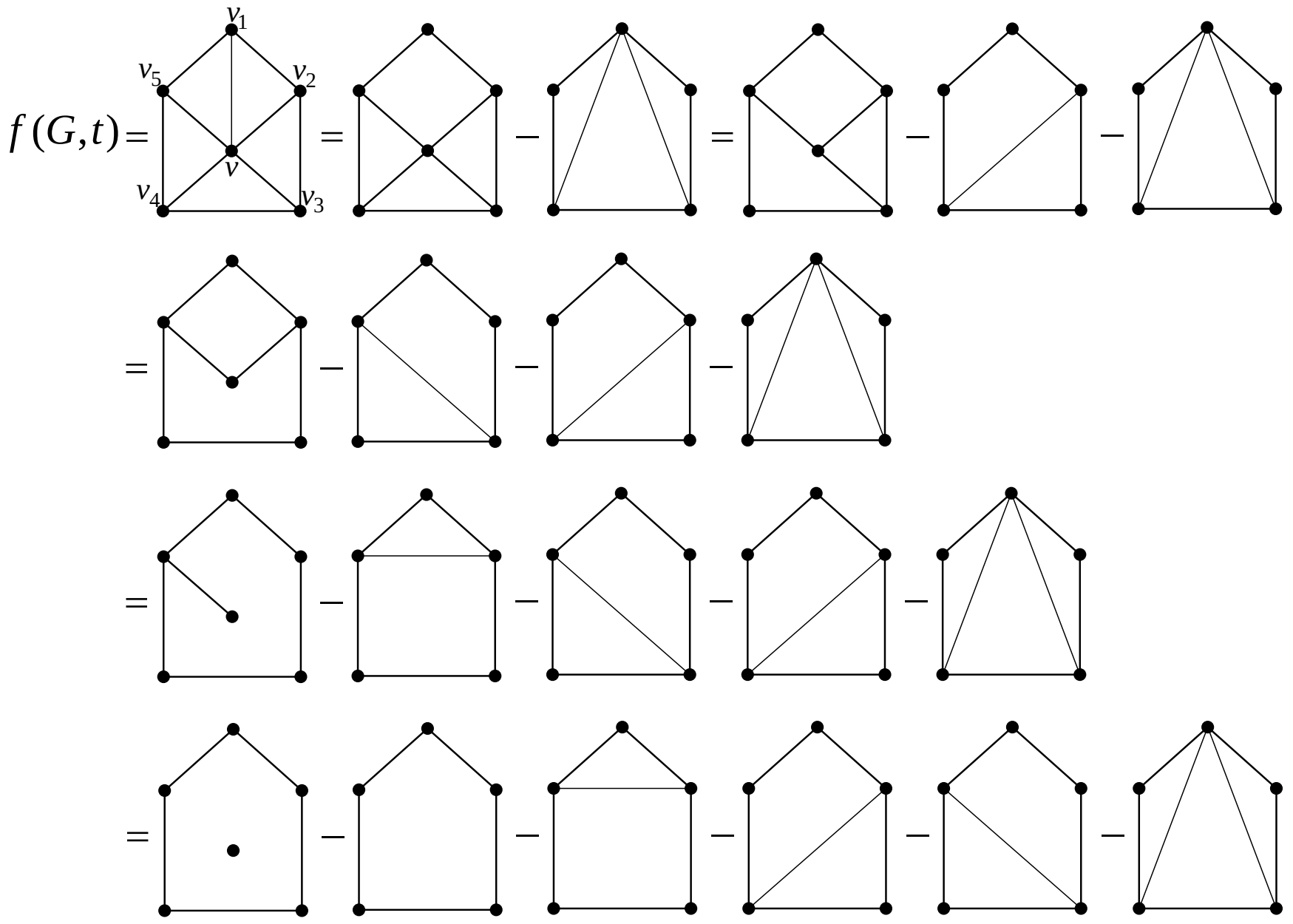

(10)

By Lemma 2, the chromatic polynomial of the first graph in Formula (10) is $t \cdot f\ G-v,t$ . Therefore, we can obtain that

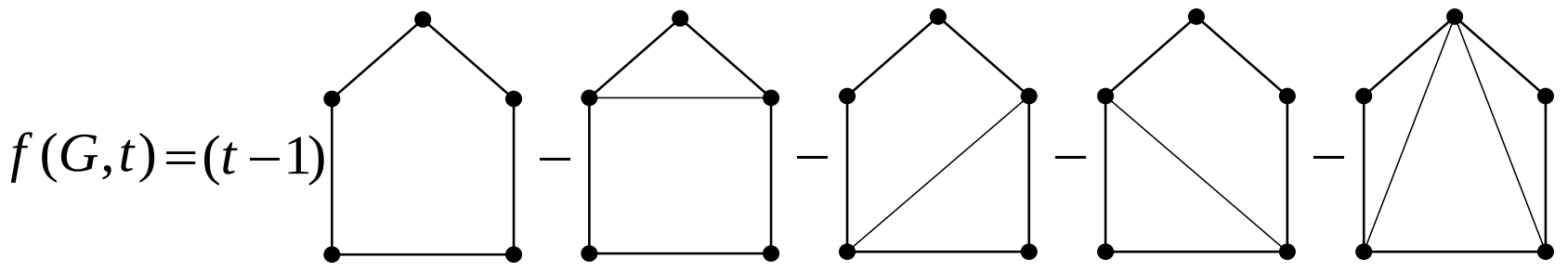

(11)

When $t=4$, we have

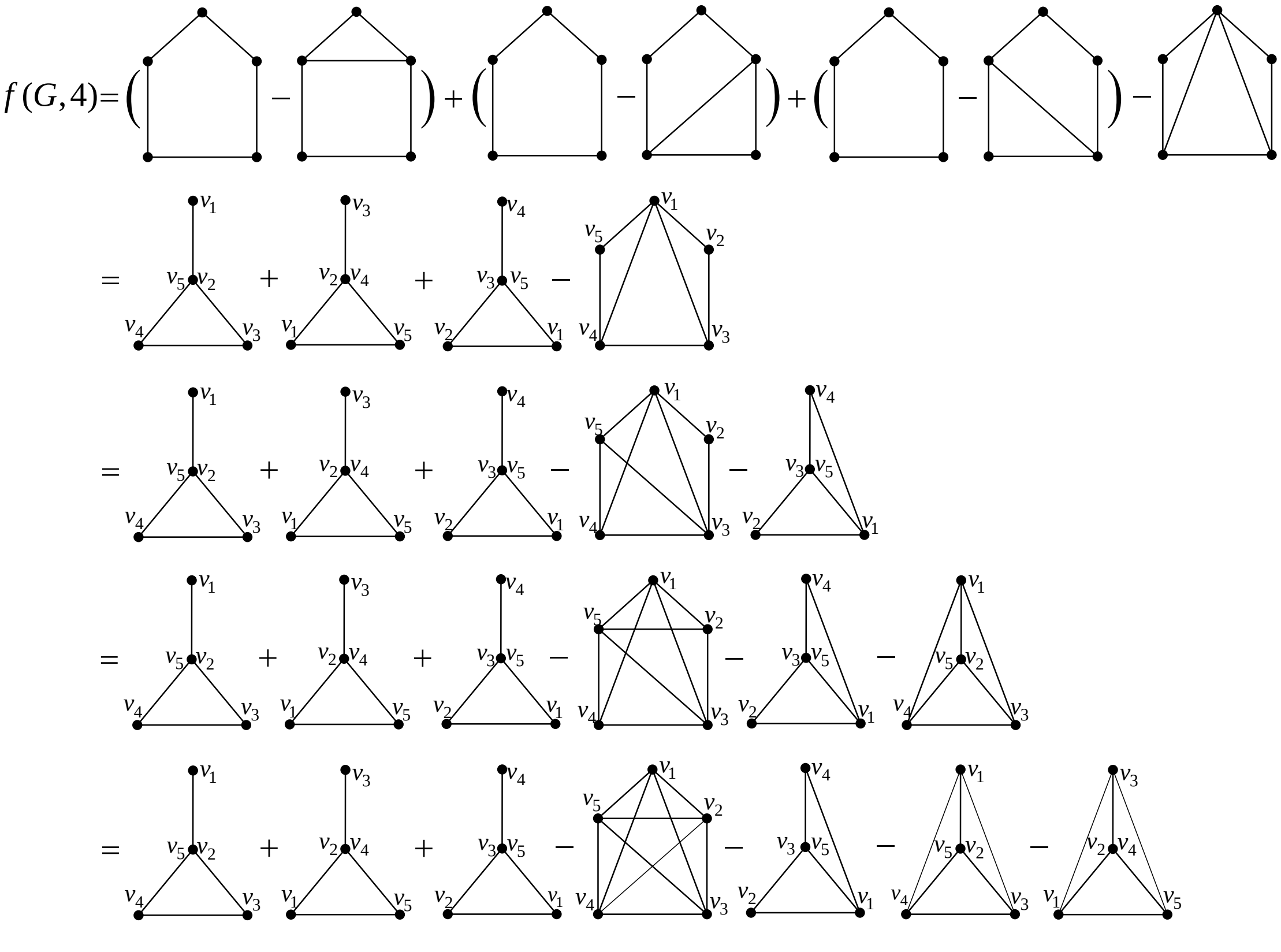

(12)

Notice that the fourth graph in Formula (12), denoted by $G'$, contains a subgraph $K_5$, and so $f(G',4)=0$. Thus, we can obtain that

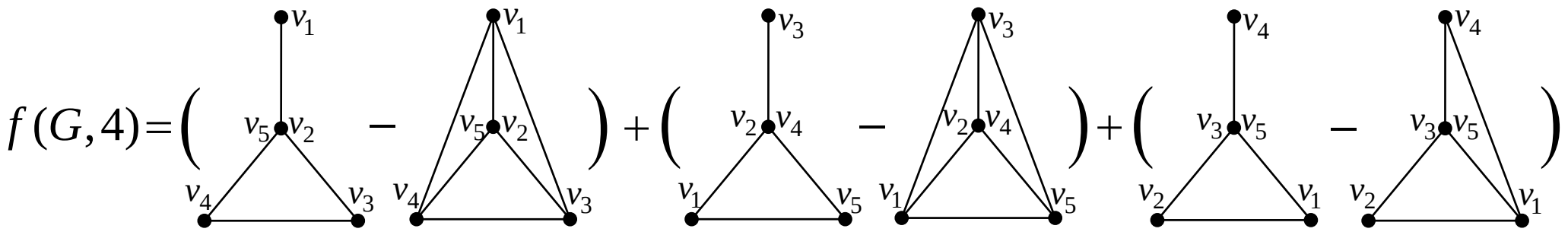

(13)

Actually, the first graph in the first bracket of Formula (13) is $G_1=(G-v)\circ\{v_2,v_5\}$; the first graph in the second bracket is $G_2=(G-v)\circ\{v_2,v_4\}$; and the first graph in the third bracket is $G_3=(G-v)\circ\{v_3,v_5\}$. The proof is completed.

## 3 Two mathematical ideas for attacking Four-Color Conjecture based on Theorem 2

It is well-known that mathematical induction is an effective method to prove Four-Color Conjecture, in which maximal planar graphs are classified into three cases by their minimum degrees. The case of minimum degree 3 or 4 is easy to prove by induction, but for the case of minimum degree 5 no mathematical method has been found. Based on Theorems 1 and 2, we give a new method to prove Four-Color Conjecture as follows.

In order to prove $f(G,4)>0$ for a maximal planar graph $G$, we use a mathematical inductive method on the number $n$ of vertices of $G$.

When $n=3,4,5$, the result is obviously true.

Assume that $n\geq 5$ and the result is true for any maximal planar graph of order at most $n-1$. We consider the case that the order of graphs is $n$. We only consider simple maximal planar graphs. Notice that for any maximal planar graph $G$, $3\leq\delta(G)\leq 5$. So we need to consider the following three cases based on the minimum degree.

**Case 1**. $\delta(G)=3$;

Let $v\in V(G)$, $d(v)=3$, and $G_1=G[\overline{N(v)}]$, $G_2=G-v$. Then we can obtain that $G_1\cap G_2=G[N(v)]\cong K_3$. Notice that $G_1=G[\overline{N(v)}]\cong K_4$, we can obtain the following result by Lemma 2:

$$f(G,t)=f(G_1\cup G_2,t)=\frac{f(G_1,t)\times f(G_2,t)}{f(K_3,t)}$$
$$=(t-3)f(G_2,t)$$

By the induction hypothesis, $f(G_2,4)>0$. Thus, $f(G,4)=f(G_2,4)>0$.

Hence, the result is true when $\delta(G)=3$.

**Case 2**. $\delta(G)=4$;

Let $v\in V(G)$, $d(v)=4$ and $N(v)=\{v_1,v_2,v_3,v_4\}$ (see Figure 2). Notice that we use $G[\overline{N(v)}]$ to denote $G$. By Theorem 1, we have $f(G,4)=f(G_1,4)+f(G_2,4)$, where $G_1=(G-v)\circ\{v_1,v_3\}$ and $G_2=(G-v)\circ\{v_2,v_4\}$. It is easy to prove that the graphs $G_1$ and $G_2$ are both maximal planar graphs with $n-2$ vertices. By the induction hypothesis, we can obtain

$$f(G_1,4)=f((G-v)\circ\{v_1,v_3\},4)>0$$
$$f(G_2,4)=f((G-v)\circ\{v_2,v_4\},4)>0$$

Therefore, $f(G,4)=f(G_1,4)+f(G_2,4)>0$ and the result is true when $\delta(G)=4$.

The key ingredient of the proof is the following Case 3.

**Case 3**. $\delta(G)=5$;

The maximal planar graph of minimum degree 5 with fewest vertices is the icosahedron, depicted in Figure 4(a), which has 12 vertices. Obviously, the icosahedron is 4-colorable. There is no maximal planar graph of minimum degree 5 with 13 vertices. Notice that for any maximal planar graph $G$ of order at least 14 and minimum degree 5, there exists a vertex $v\in V(G)$ such that $d(v)=5$ and $d(v_1)\geq 6$, where $N(v)=\{v_1,v_2,v_3,v_4,v_5\}$ (see Figure 3). Hence, the graph $G_1$ in Theorem 2 is a 4-colorable maximal planar graph of minimum degree at least 4. Based on this evidence, we give two mathematical ideas to prove $f(G,4)>0$ as follows.

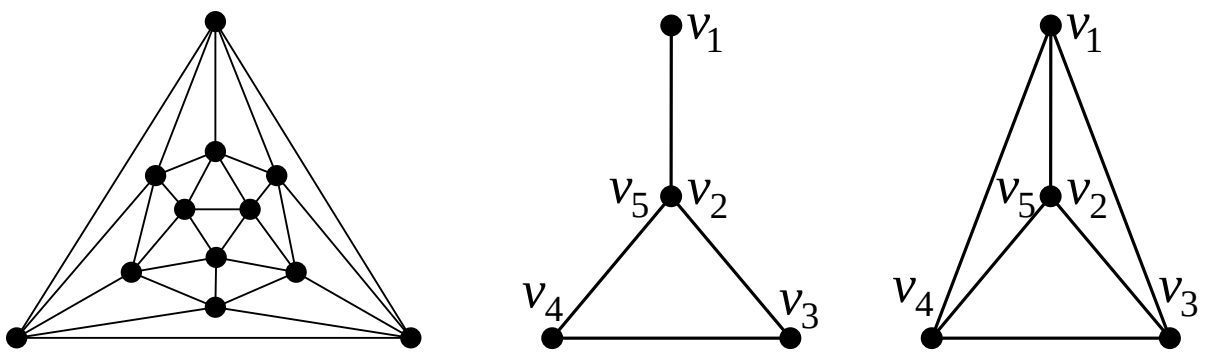

(a) The icosahedron (b) A funnel subgragh (c)

Figure 4. Three graphs in Case 3.

The first idea is based on the fact that the value of each square bracket in Formula (9) is no less than zero. Hence, the Four-Color Conjecture can be proved if any square bracket's value is greater than zero. The value of the first square bracket is greater than zero if and only if there exists $f_1 \in C_4^0\ G_1$ such that $f_1\ v_1 = f_1\ v_3$ or $f_1\ v_1 = f_1\ v_4$. Therefore, $f(G,4)=0$ if and only if each square bracket in Formula (9) is equal to zero. Moreover, the value of the first square bracket is zero if and only if for any $f_1 \in C_4^0\ G_1$, $f_1\ v_1 \neq f_1\ v_3$ and $f_1\ v_1 \neq f_1\ v_4$, that is, for any $f_1 \in C_4^0\ G_1$, the colors of vertices of the funnel shown in Figure 4(b) are pairwise different. Such maximal planar graphs are called *4-chromatic-funnel pseudo uniquely-4-colorable maximal planar graphs*. For instance, each graph in Figure 5 is a 4-chromatic-funnel pseudo-uniquely 4-colorable maximal planar graph.

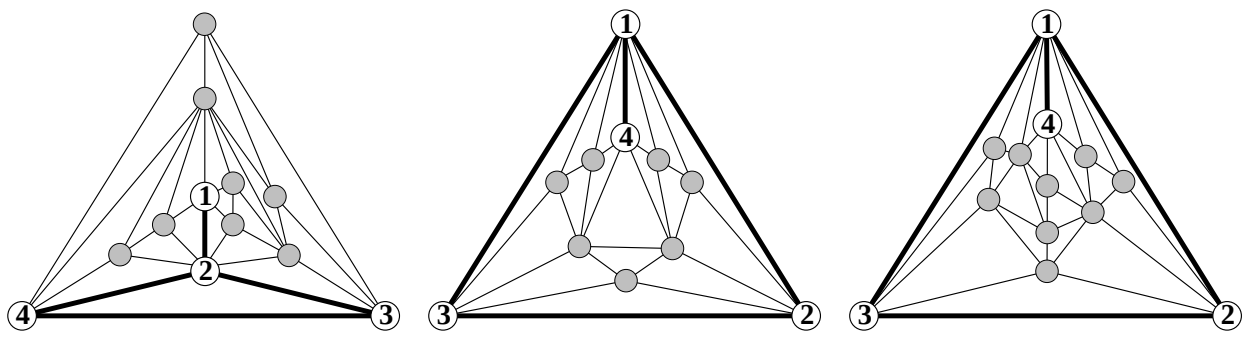

(a) A 4-chromatic-funnel type graph of order 12 (b) A 4-chromatic-funnel type graph of order 11 (c) A 4-chromatic-funnel type graph of order 13

Figure 5. Three 4-chromatic-funnel pseudo uniquely-4-colorable maximal planar graphs.

A $k$-colorable graph $G$ is called a *$k$-colorable coordinated graph* if there exist $k$ vertices $v_1, v_2, \cdots, v_k$ in $G$ such that $f\ v_1\ , f\ v_2\ , \cdots, f\ v_k$ are pairwise different for any $k$-coloring $f$ of $G$. All 4-colorable coordinated maximal planar graphs can be divided into three classes: (1) uniquely 4-colorable maximal planar graphs, namely these graphs have only one partition of $k$-color class; (2) quasi uniquely-4-colorable maximal planar graphs, namely these graphs contain a subgraph that is uniquely 4-colorable; (3) pseudo uniquely-4-colorable maximal planar graphs, namely these graphs that are neither uniquely 4-colorable nor quasi uniquely-4-colorable. A detailed research on 4-colorable coordinated maximal planar graphs will be given in the subsequent series of articles.

Now we give the second idea to prove $f\ G,4\ > 0$. The maximal planar graphs $G_1, G_2$ and $G_3$ in Theorem 2 can be regarded as the graphs obtained from $G$ by deleting a 5-degree vertex $v$ and contracting $v_2, v_5$, $v_2, v_4$ and $v_3, v_5$ into a single vertex, respectively (see Figure 6). Moreover, the 5-cycle consisting of the neighbors of $v$ in $G$ is contracted to a funnel subgraph $L_1 = v_1 - \Delta v_2^5 v_3 v_4$ in $G_1$, $L_2 = v_3 - \Delta v_2^4 v_1 v_5$ in $G_2$, and $L_3 = v_4 - \Delta v_3^5 v_1 v_2$ in $G_3$, respectively, where $v_2^5$, $v_2^4$ and $v_3^5$ are the new vertices obtained by contracting $v_2, v_5$, $v_2, v_4$ and $v_3, v_5$, respectively.

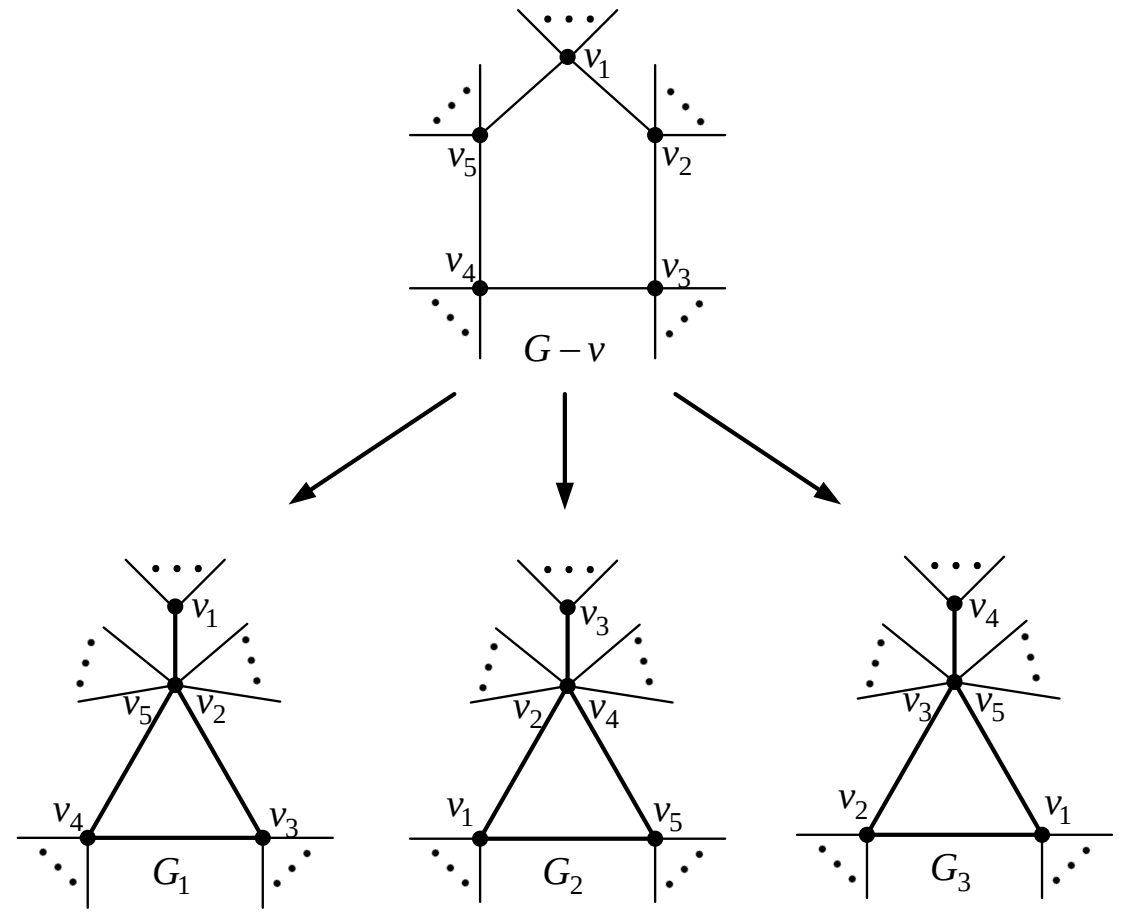

Figure 6. The processes of generating the three funnel subgraphs

By the induction hypothesis, $G_1$, $G_2$ and $G_3$ are 4-colorable. In order to prove $f\ G,4\ > 0$, it is needed to prove that at least one of the funnel subgraphs $L_1$, $L_2$ and $L_3$ is not 4-chromatic.

Therefore, the second idea is to prove that for any maximal planar graph $G$ of minimum degree 5, there exists a 5-wheel $W_5^v$ in $G$ such that at least one of the funnel subgraphs $L_1$, $L_2$ and $L_3$ corresponding to $G_1, G_2$ and $G_3$ is not a 4-chromatic-funnel. For instance, the graph in Figure 5(a) can be regarded as the maximal planar graph obtained from the graph in Figure 7 by the operation shown in Figure 6. It is not difficult to prove that the other two graphs obtained from Figure 7 by the operation shown in Figure 6 have no 4-chromatic-funnel.

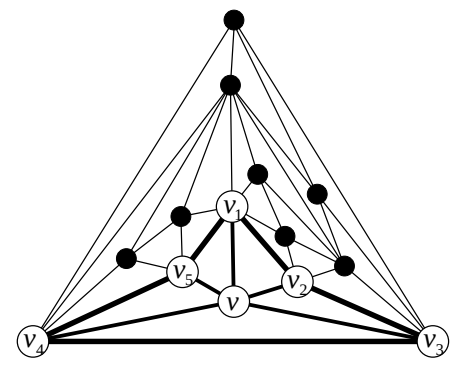

Figure 7. A maximal planar graph that can be contracted to the graph in Figure 5(a).

## 4 Conclusion

In this paper we give two recursion formulae of chromatic polynomial on maximal planar graphs. Based on these formulae, we find: (1) two mathematical ideas for attacking Four-Color Conjecture; (2) a method to generate maximal planar graphs, called contracting and extending operational system, which establishes a relation between the structure and colorings of a maximal planar graph. For instance, the maximal planar graph in Figure 5(a) can be obtained from the graph in Figure 7 by the extending 5-wheel operation, in other words, the maximal planar graph in Figure 7 can be obtained from the graph in Figure 5(a) by the contracting 5-wheel operation. A detailed research on contracting and extending operational system of maximal planar graphs will be given in later articles.

**Acknowledgements**: Thanks to my students Enqiang Zhu, Zepeng Li, Xiaoqing Liu, Hongyu Wang and Yangyang Zhou for some useful discussions.

References

[1] JENSEN T R and TOFT B. Graph Coloring Problems[M]. New York: John Wiley & Sons, 1995: 48–49.

[2] JOSEP DÍAZ, JORDI PETIT, and MARIA SERNA. A survey of graph layout problems[J]. *ACM Computing Surveys*, 2002, 34(3): 313–355.

[3] BRODER A, KUMAR R, MAGHOUL F, *et al.* Graph structure in the Web[J]. *Computer Networks*, 2000, 33(1–6): 309–320.

[4] XU J, LI Z P and ZHU E Q. Research progress on the theory of maximal planar graphs[J]. *Chinese Journal of Computers*, 2015, 38(7): 1680–1704. (in Chinese)

[5] KEMPE A B. On the geographical problem of the four colors [J]. *American Journal of Mathematics*, 1879, 2(3): 193–200.

[6] HEAWOOD P J. Map colour theorem[J]. *Quarterly Journal of Mathematics*, 1890, 24: 332–338.

[7] APPEL K and HAKEN W. The solution of the four-color map problem[J]. *Science Ameri*can, 1977, 237(4): 108–121.

[8] APPEL K and HAKEN W. Every planar map is four colorable, I: Discharging[J]. *Illinois Journal of Mathematics*, 1977, 21(3): 429–490.

[9] APPEL K, HAKEN W and KOCH J. Every planar map is four-colorable, II: reducibility[J]. *Illinois Journal of Mathematics*, 1977, 21(3): 491–567.

[10] ROBERTSON N, SANDERS D P, SEYMOUR P, *et al.* A new proof of the four colour theorem[J]. *Electronic Research Announcements American Mathematical Society*, 1996, 2: 17–25.

[11] ROBERTSON N, SANDERS D P, SEYMOUR P D, *et al.* The four color theorem[J]. *Journal of Combinatorial Theory*, Series B, 1997, 70(1): 2–44.

[12] WERNICKE P. Über den kartographischen Vierfarbensatz [J]. *Mathematische Annalen*, 1904, 58(3): 413–426.

[13] BIRKHOFF G D. The reducibility of maps[J]. *American Journal of Mathematics*, 1913, 35(2): 115–128.

[14] HEESCH H. Untersuchungen Zum Vierfarbenproblem[M]. Mannheim/Wien/Zürich: Bibliographisches Institut, 1969: 4–12.

[15] FRANKLIN P. The four color problem[J]. *American Journal of Mathematics*, 1922, 44(3): 225–236.

[16] FRANKLIN P. Note on the four color problem[J]. *Journal of Mathematical Physics*, 1938, 16: 172–184.

[17] REYNOLDS C. On the problem of coloring maps in four colors[J]. *Annals of Mathematics*, 1926–27, 28(1–4): 477–492.

[18] WINN C E. On the minimum number of polygons in an irreducible map[J]. *American Journal of Mathematics*, 1940, 62(1): 406–416.

[19] ORE O and STEMPLE J. Numerical calculations on the four-color problem[J]. *Journal of Combinatorial Theory*, 1970, 8(1): 65–78.

[20] MAYER J. Une propriété des graphes minimaux dans le problème des quatre couleurs[J]. *Problèmes Combinatoires et Thorie des Graphes*, *Colloques Internationaux CNRS*, 1978, 260: 291–295.

[21] TAIT P G. Remarks on the colouring of maps[J]. *Proceedings of the Royal Society of Edinburgh*, 1880, 10: 501–516.

[22] PETERSEN J. Surle théorème de Tait[J]. *L'intermédiaire des Mathématiciens*, 1898, 5: 225–227.

[23] TUTTE W T. On Hamiltonian circuits[J]. *Journal London Mathematical Society*, 1946, 21: 98–101.

[24] GRINBERG E J. Plane homogeneous graphs of degree three without Hamiltonian circuits[J]. *Latvian Math Yearbook*, 1968, 5: 51–58.

[25] BIRKHOFF G D. A determinantal formula for the number of ways of coloring a map[J]. *Annals of Mathematics*, 1912, 14: 42–46.

[26] BIRKHOFF G D and LEWIS D. Chromatic polynomials[J]. *Transactions of the American Mathematical Society*, 1946, 60: 355–451.

[27] DONG F M, KOH K M, and TEO K L. Chromatic Polynomials and Chromaticity of Graphs[M]. World Scientific, Singapore, 2005: 23–215.

[28] TUTTE W T. On chromatic polynomials and the golden ratio[J]. *Journal of Combinatorial Theory*, 1970, 9(3):

289–296.

[29] TUTTE W T. Chromatic sums for planar triangulations, V: Special equations[J]. *Canadian Journal of Mathematics*, 1974, 26: 893–907.

[30] READ R C. An introduction to chromatic polynomials[J]. *Journal of Combinatorial Theory*, 1968, 4(1): 52–71.

[31] WHITNEY H. On the coloring of graphs[J]. *Annals of Mathematics*, 1932, 33(4): 688–718.

[32] XU J. Recursive formula for calculating the chromatic polynomial of a graph by vertex deletion[J]. *Acta Mathematica Scientia*, 2004, 24$B$(4): 577–582.

[33] XU J and LIU Z. The chromatic polynomial between graph and its complement — About Akiyama and Hararys'open problem[J]. *Graph and Combinatorics*, 1995, 11: 337–345.

[34] ZYKOV A A. On some properties of linear complexes[J]. *Math Ussr Sbornik*, 1949, 24(66): 163–188 (in Russian); *English Translation in Transactions of the American Mathematical Society*, 1952, 79.

XU Jin： Born in 1959, Professor. His main research interests include graph theory and combinatorial optimization, biocomputing, social networks and information security.